\documentclass[reqno,12pt]{amsproc}
  \usepackage{latexsym} 
  \usepackage[all]{xy}
  \usepackage{amsfonts} 
  \usepackage{amsthm} 
  \usepackage{amsmath} 
  \usepackage{amssymb}
  \usepackage{pifont}  
  
  \usepackage[T1]{fontenc}

  \def\sw#1{{\sb{(#1)}}} 
   
  \def\su#1{{\sp{[#1]}}}  
  \def\suc#1{{\sp{(#1)}}}

  \def\<{{\langle}} 
  \def\>{{\rangle}} 
  \def\ra{{\triangleleft}} 
   
  \def\eps{\varepsilon}

  \def\note#1{{}}

\def\flip{{\rm \sf{\; flip\; }}}

  \def\note#1{}

  \def\cM{{\mathfrak M}}

  \def\beq{\begin{equation}} 
  \def\eeq{\end{equation}}

  \def\id{\mathrm{id}}

  \def\ot{{\otimes}}

  \def\roM{\varrho^{M}} 
    \def\rom{\varrho_{M}} 
 \def\lam{\widehat{\varrho}_{M}} 
  \def\lan{\widehat{\varrho}_{N}} 
   \def\laM{\widehat{\varrho}^{M}} 
  \def\laN{\widehat{\varrho}^{N}}

    \def\vect{\mathfrak{Vect}}
     \def\1{\mathbf{1}}

\def\uH{\underline{H}}
\def\oH{\overline{H}}
\def\act{\!\cdot\!}

\def\k{\Bbbk}

  \newcounter{zlist} 
  \newenvironment{zlist}{\begin{list}{(\arabic{zlist})}{ 
  \usecounter{zlist}\leftmargin2.5em\labelwidth2em\labelsep0.5em 
  \topsep0.6ex%\itemsep0.3ex plus0.2ex minus0.3ex 
  \parsep0.3ex plus0.2ex minus0.1ex}}{\end{list}}

  \newcounter{blist} 
  \newenvironment{blist}{\begin{list}{(\alph{blist})}{ 
  \usecounter{blist}\leftmargin2.5em\labelwidth2em\labelsep0.5em 
  \topsep0.6ex %\itemsep0.3ex plus0.2ex minus0.3ex 
  \parsep0.3ex plus0.2ex minus0.1ex}}{\end{list}} 

  \newcounter{rlist}

   \newcounter{alist}

\def\stac#1{\raise-.2cm\hbox{$\stackrel{\displaystyle\otimes}{\scriptscriptstyle{#1}}$}}

\def\cten#1{\raise-.2cm\hbox{$\stackrel{\displaystyle\widehat{\otimes}}
{\scriptscriptstyle{#1}}$}}

  \def\Label#1{\label{#1}\ifmmode\llap{[#1] }\else 
  \marginpar{\smash{\hbox{\tiny [#1]}}}\fi} 
  \def\Label{\label}

  \newtheorem{proposition}{Proposition}[section]
  \newtheorem{lemma}[proposition]{Lemma}
  \newtheorem{corollary}[proposition]{Corollary} 
  \newtheorem{theorem}[proposition]{Theorem}

  \theoremstyle{definition} 
  \newtheorem{definition}[proposition]{Definition}

  \theoremstyle{remark}

  \theoremstyle{definition} 
  \swapnumbers

\begin{document} 

 \title{Hopf modules and the fundamental theorem for Hopf (co)quasigroups} 
 \author{Tomasz Brzezi\'nski}
 \address{ Department of Mathematics, Swansea University, 
  Singleton Park, 
 \newline\indent  
  Swansea SA2 8PP, U.K.} 
  \email{T.Brzezinski@swansea.ac.uk}   
    \date{December 2009} 
  \begin{abstract} 
The notion of a Hopf module over a Hopf (co)quasigroup is introduced and a version of the fundamental theorem for Hopf (co)quasigroups is proven.
  \end{abstract} 
  \maketitle

\section{Introduction}
Generalisations of Hopf algebras over fields have quite a long history. The main common principle of  many (but not all) such generalisations is the weakening of some of the algebraic conditions that enter the definition of a Hopf algebra. We mention here but a few examples. The weakening of the (co)associativity leads to {\em quasi-Hopf algebras} \cite{Dri:qua}. If the algebra part of a Hopf algebra is not required to have a unit, one obtains {\em multiplier Hopf algebras} \cite{Van:mul}.  If one does not require the unit to be comultiplicative and the counit to be multiplicative, one is led to {\em weak Hopf algebras} \cite{BohNil:wea}. The most recent additions to this family of generalisations are {\em Hopf quasigroups} and {\em Hopf coquasigroups} introduced in \cite{KliMaj:Hop} in order to capture the quasigroup features of the (algebraic) 7-sphere. Similarly to quasi-Hopf algebras, Hopf (co)quasigroups are not required to be (co)associative. The lack of (co)associativity is compensated by conditions involving the antipode. 

The first aim of this note is to show that, similarly to standard Hopf algebras, Hopf (co)\-quasi\-groups can be characterised by a Galois-type condition. More specifically, that  a Hopf coquasigroup satisfies the right Galois condition, i.e.\ that the right Galois map (see Definition~\ref{def.Gal}) is bijective,  is already proven in \cite[Lemma~6.3]{KliMaj:Hop}. Parallel arguments yield satisfaction of the left Galois condition for Hopf coquasigroups and both Galois conditions for Hopf quasigroups. We provide the converse to this statement, namely that one needs both left and right Galois maps to have almost (co)linear inverses (see Definition~\ref{def.almost}) to infer the existence of an antipode for a Hopf (co)quasigroup. We also lay foundations for the theory of Hopf modules over Hopf (co)quasigroups, and show in particular that the categories of such modules are equivalent to the category of vector spaces. All these results can be understood as a Hopf (co)quasigroup version of the fundamental theorem of Hopf algebra theory. They can also be interpreted as the general background that allows for the development of differential structures on Hopf coquasigroups presented in \cite[Section~6]{KliMaj:Hop}. 

All algebras and coalgebras are over a field $\k$. Unadorned tensor product symbol represents the tensor product of $\k$-vector spaces. The identity map on a vector space $V$ is denoted by $V$. The unit of a Hopf algebra $H$ both as an element of $H$ and as a map $\k\to H$ is denoted by 1. The product of elements of a Hopf algebra is denoted by juxtaposition. 

\section{The first fundamental theorem for Hopf (co)quasigroups}\label{sec.fund}
\setcounter{equation}{0}
The aim of this section is to show that the definition of a Hopf (co)quasigroup can be rephrased in terms of bijectivity of (canonical) Galois maps.

\begin{definition}[\cite{KliMaj:Hop}]\label{def.quasi}
Let $H$ be a vector space that is a unital (not necessarily associative) algebra with product $\mu: H\ot H\to H$ and unit $1: \k \to H$, and a counital (not necessarily coassociative) coalgebra with coproduct $\Delta: H\to H\ot H$ and counit $\eps: H\to \k$ that are algebra homomorphisms. 

$H$ is called a {\em Hopf quasigroup} provided $\Delta$ is coassociative and there exists a linear map $S: H\to H$ such that
\begin{equation}\label{quasi1}
\mu\circ (H \ot \mu) \circ (S\ot H\ot H)\circ (\Delta \ot H) = \eps \ot H = \mu\circ (H \ot \mu) \circ (H\ot S\ot H)\circ (\Delta \ot H)
\end{equation}
and
\begin{equation} \label{quasi2}
\mu\circ (\mu \ot H) \circ (H\ot H\ot S)\circ (H\ot \Delta) = H\ot \eps = \mu\circ (\mu \ot H) \circ (H\ot S\ot H)\circ (H\ot \Delta).
\end{equation}

$H$ is called a {\em Hopf coquasigroup} provided $\mu$ is associative and there exists a linear map $S: H\to H$ such that
\begin{equation}\label{coq1}
(\mu \ot H) \circ (S\ot H\ot H)\circ (H\ot \Delta)\circ \Delta = 1\ot H = (\mu \ot H) \circ (H\ot S\ot H)\circ (H\ot \Delta)\circ \Delta
\end{equation}
and
\begin{equation}\label{coq2}
(H \ot \mu) \circ (H\ot H\ot S)\circ (\Delta \ot H)\circ \Delta = H\ot 1 = (H \ot \mu) \circ (H\ot S\ot H)\circ (\Delta \ot H)\circ \Delta .
\end{equation}
\end{definition}

We use Sweedler notation for coproduct, that is, for all $h\in H$, $\Delta (h) = h\sw 1\ot h\sw 2$. One should remember, however, that if $\Delta$ is not coassociative, the standard Sweedler's relabeling rules no longer apply.  As for standard Hopf algebras, the map $S$ in Definition~\ref{def.quasi} is called an {\em antipode}. To relieve the notation from  some brackets, we write $Sh$ for the value of $S$ at $h$. It is proven in \cite{KliMaj:Hop} that the antipode is antimultiplicative and anticomultiplicative and it immediately follows from (any of) equations \eqref{quasi1}--\eqref{coq2} that, for all $h\in H$,
$
(Sh\sw 1)h\sw 2 = h\sw 1 Sh\sw 2 = \eps(h)1,
$
i.e.\ $S$ enjoys the standard antipode property. 

\begin{definition}\label{def.almost}
Let $H$ be a (not necessarily) associative algebra with a compatible (not necessarily coassociative) coalgebra structure as in Definition~\ref{def.quasi}. Consider a $\k$-linear map $\phi: H\ot H\to H\ot H$. We say that
\begin{blist}
\item $\phi$ is {\em almost left $H$-linear} if, for all $g,h\in H$,
$
\phi(g\ot h) = (g\ot 1)\phi(1\ot h);
$
\item $\phi$ is {\em almost right $H$-linear} if, for all $g,h\in H$,
$
\phi(g\ot h) = \phi(g\ot 1)(1\ot h);
$
\item $\phi$ is {\em almost left $H$-colinear} if
$
\phi = (H\ot \eps\ot H)\circ (H\ot \phi)\circ (\Delta\ot H),
$
that is, for all $g,h\in H$, $\phi(g\ot h) = g\sw 1 \ot (\eps \ot H)(\phi(g\sw 2\ot h))$;
\item $\phi$ is {\em almost right $H$-colinear} if
$
\phi = (H\ot \eps\ot H)\circ (\phi \ot H)\circ (H\ot \Delta),
$
that is, for all $g,h\in H$, $\phi(g\ot h) = (H \ot \eps)(\phi(g\ot h\sw 1))\ot h\sw 2$.
\end{blist}
\end{definition}
In case (co)product of $H$ is (co)associative, the notion of almost $H$-(co)linearity coincides with that of  $H$-(co)linearity. 

\begin{definition}\label{def.Gal}
Let $H$ be a (not necessarily) associative algebra with a compatible (not necessarily coassociative) coalgebra structure as in Definition~\ref{def.quasi}. The $\k$-linear map
\begin{equation}\label{beta}
\beta : H\ot H \to H\ot H, \qquad \beta= (\mu\ot H)\circ (H\ot \Delta) : g\ot h \mapsto gh\sw 1\ot h\sw 2,
\end{equation}
is called the {\em right Galois map}. The $\k$-linear map
\begin{equation}\label{gamma}
\gamma : H\ot H \to H\ot H, \qquad \gamma= (H\ot \mu)\circ (\Delta\ot H) : g\ot h \mapsto g\sw 1\ot g\sw 2h,
\end{equation}
is called the {\em left Galois map}. 
\end{definition}

\begin{lemma}\label{lemma.almost}
The right Galois map $\beta$ is almost left $H$-linear and almost right $H$-colinear. The left Galois map $\gamma$ is almost right $H$-linear and almost left $H$-colinear.
\end{lemma}
\begin{proof}
The statements are obviously related by the left-right symmetry, hence we prove only the first one. The fact that $\beta$ is almost left $H$-linear follows immediately from the observation that, for all $h\in H$, $\beta(1\ot h) = \Delta(h)$, and from the unitality of $\mu$. In view of the counitality of $\Delta$, we obtain, for all $g,h \in H$,
$$
(H \ot \eps)(\beta(g\ot h\sw 1))\ot h\sw 2 = gh\sw 1\sw 1\ot \eps(h\sw 1\sw 2) \ot h\sw 2 = gh\sw 1\ot h\sw 2 = \beta(g\ot h),
$$
hence $\beta$ is almost right $H$-colinear as required. 
\end{proof}

The main result of this section is contained in the following (one direction, that $\beta$ is bijective for a Hopf coquasigroup is already proven in \cite[Lemma~6.3]{KliMaj:Hop})

\begin{theorem}\label{thm.fund}
Let $H$ be a vector space that is a unital (not necessarily associative) algebra with product $\mu: H\ot H\to H$ and unit $1: \k \to H$, and a counital (not necessarily coassociative) coalgebra with coproduct $\Delta: H\to H\ot H$ and counit $\eps: H\to \k$ that are algebra homomorphisms. 
\begin{zlist}
\item $H$ is a Hopf quasigroup if and only if $\Delta$ is coassociative and the right and left Galois maps have almost left, resp.\ right $H$-linear inverses.
\item $H$ is a Hopf coquasigroup if and only if $\mu$ is associative and the right and left Galois maps have almost right, resp.\ left $H$-colinear inverses.
\end{zlist}
\end{theorem}
\begin{proof}
(1) If $H$ is a Hopf quasigroup, then, by definition, $\Delta$ is coassociative. Set
\begin{equation}\label{beta-1}
\beta^{-1}: H\ot H\to H\ot H, \qquad g\ot h\mapsto gSh\sw 1 \ot h\sw 2,
\end{equation}
and
\begin{equation}\label{gamma-1}
\gamma^{-1}: H\ot H\to H\ot H, \qquad g\ot h\mapsto g\sw 1 \ot (Sg\sw 2)h.
\end{equation}
The map $\beta^{-1}$ is almost left $H$-linear and $\gamma^{-1}$ is almost right $H$-linear. Furthermore, that $\beta^{-1}$ is the inverse to $\beta$ follows by a calculation parallel to that in the proof of \cite[Lemma~6.3]{KliMaj:Hop}, which establishes the bijectivity of $\beta$ for a Hopf coquasigroup. Explicitly, for all $g,h\in H$,
$$
\beta^{-1}\circ \beta (g\ot h) =\beta^{-1}(gh\sw 1\ot h\sw 2) = (gh\sw1)Sh\sw 2 \ot h\sw 3 = g\ot h,
$$
by the first of equations \eqref{quasi2}. Similarly, using the second of equations \eqref{quasi2} one computes
$$
\beta\circ\beta^{-1}(g\ot h) = (gSh\sw 1)h\sw 2 \ot h\sw 3 = g\ot h.
$$
Therefore, $\beta^{-1}$ is the inverse of $\beta$. Analogous computation that uses equations \eqref{quasi1} yields that $\gamma^{-1}$ is the inverse of $\gamma$.

Conversely, assume that the Galois maps $\beta$ and $\gamma$ have the required inverses and introduce the notation, for all $h\in H$,
$$
h\su 1 \ot  h\su 2 :=  \beta^{-1}(1\ot h), \qquad h\suc 1 \ot h\suc 2 := \gamma^{-1}(h\ot 1).
$$
Since $\beta^{-1}$ and $\gamma^{-1}$ are almost $H$-linear, for all $g,h\in H$,
\begin{equation}\label{trans}
\beta^{-1} (g\ot h) = gh\su 1\ot h\su 2, \qquad \gamma^{-1}(h\ot g) = h\suc 1\ot h\suc 2g.
\end{equation}
Define $\k$-linear maps $S: H\to H$, $\bar{S}: H\to H$, by
\begin{equation}\label{s}
S: h\mapsto h\su 1\eps(h\su 2), \qquad \bar{S} : h\mapsto \eps(h\suc 1)h\suc 2,
\end{equation}
that is 
\begin{equation}\label{scomp}
S = (H\ot \eps)\circ \beta^{-1} \circ (1\ot H), \qquad \bar{S} = (\eps\ot H)\circ \gamma^{-1}\circ  (H\ot 1).
\end{equation}
Since the coproduct $\Delta$ is coassociative, the almost right colinearity of $\beta$ is tantamount with its colinearity, and hence also $\beta^{-1}$ is right $H$-colinear. This implies that, for all $h\in H$,
$$
h\sw 1\su 1 \ot h\sw 1\su 2 \ot h\sw 2 = h\su 1\ot h\su 2\sw 1\ot h\su 2\sw 2.
$$
Applying the counit to the middle factor one concludes that
$$
h\su 1\ot h\su 2 = h\sw 1\su 1\eps(h\sw 1\su 2)\ot h\sw 2 = Sh\sw 1\ot h\sw 2.
$$
In view of equation \eqref{trans}, the map $\beta^{-1}$ has the same form as in \eqref{beta-1}. Calculations analogous to those in the first part (i.e.\ the use of the fact that $\beta^{-1}$ is the inverse of $\beta$) produce the equalities
$$
(gh\sw1)Sh\sw 2 \ot h\sw 3 = g\ot h = (gSh\sw 1)h\sw 2 \ot h\sw 3.
$$
Applying the counit to the second factor, one obtains equations \eqref{quasi2}. 

Following similar chain of arguments one concludes that the map $\bar{S}$ satisfies equations \eqref{quasi1}. Finally, using the facts that $S$ satisfies equations \eqref{quasi2} (in the second equality), and that $\bar{S}$ satisfies equalities \eqref{quasi1} (in the third step), we can compute, for all $h\in H$,
$$
\bar{S}(h) = \bar{S}h\sw 1\eps(h\sw 2) = ((\bar{S}h\sw 1)h\sw 2)Sh\sw 3 = \eps(h\sw 1)Sh\sw 2 = Sh.
$$
Thus the maps $S$ and $\bar{S}$ coincide and they are the required antipode of the Hopf quasigroup. 

(2) Note that the definition of $S$ and $\bar{S}$  in equations \eqref{scomp} is self-dual under the usual duality (reversing arrows, interchanging algebra with coalgebra structures), which interchanges Hopf quasigroups with Hopf coquasigroups. Thus the similar reasoning as in part (1) gives the proof of part (2).
\end{proof}

\section{Hopf modules over Hopf (co)quasigroups}\label{sec.mod} \setcounter{equation}{0}
The aim of this section is to introduce Hopf modules over Hopf (co)quasigroups and to show that, similarly to ordinary Hopf algebras, the categories of such modules are equivalent to the category of vector spaces.

\subsection{Hopf modules over Hopf quasigroups}
\begin{definition}\label{def.mod}
Let $H$ be a Hopf quasigroup. Let $M$ be a coassociative and counital right $H$-comodule with coaction $\roM$ and a unital right $H$-module with action $\varrho_M$. $M$ is said to be a {\em right $H$-Hopf module} if
\begin{equation}\label{mod1}
\rom\circ (\rom \ot H) \circ (M\ot H\ot S)\circ (M\ot \Delta) = M\ot \eps = \rom\circ (\rom \ot H)\circ (M\ot S\ot H)\circ (M\ot \Delta)
\end{equation}
and
\begin{equation}\label{mod2}
\roM\circ \rom = (\rom \ot \mu)\circ (M\ot \flip \ot H)\circ (\roM\ot \Delta).
\end{equation}
\end{definition}

Writing $m\act h$ for $\rom(m\ot h)$ and using the Sweedler notation for coaction $\roM(m) = m\sw 0 \ot m\sw 1$, the conditions \eqref{mod1} and \eqref{mod2} read explicitly, for all $h\in H$ and $m\in M$,
$$
(m\act h\sw 1)\act Sh\sw 2 = m\eps(h) = (m\act Sh\sw 1)\act h\sw 2 \quad \!\!\!\mbox{and} \quad \!\!\!
(m\act h)\sw 0 \ot (m\act h)\sw 1 = m\sw 0\act h\sw 1\ot m\sw 1h\sw 2.
$$
Equation \eqref{mod2} is the usual compatibility condition for Hopf modules \cite[Section~4.1]{Swe:Hop}, while equation \eqref{mod1} is a substitute for the associativity of the $H$-action, as it is  satisfied automatically whenever right $H$-action is associative. Note that $H$ is a right $H$-Hopf module with the multiplication as action and coproduct as coaction. 

We begin the analysis of Hopf modules over a Hopf quasigroup by the following simple (and completely classical)

\begin{lemma}\label{lem.rom1}
Let $H$ be a Hopf quasigroup and let $M$ be a right $H$-Hopf module. Then,
for all $m\in M$,
\begin{equation}\label{rom1}
\roM(m\sw 0 \act Sm\sw 1) = m\sw 0 \act Sm\sw 1 \ot 1.
\end{equation}
Consequently, for all $m\in M$ and $h\in H$,
\begin{equation}\label{rom2}
\roM \left( \left (m\sw 0\act Sm\sw 1\right)\act h\right) = \left (m\sw 0\act Sm\sw 1\right)\act h\sw 1 \ot h\sw 2.
\end{equation}
\end{lemma}
\begin{proof}
Equation \eqref{rom1} follows by the compatibility condition \eqref{mod2}, and by the anticomultiplicativity and the antipode property of $S$. Then equation \eqref{rom2} is a simple consequence of \eqref{rom1} and \eqref{mod2} (all as in the standard Hopf algebra case).
\end{proof}

On a Hopf module over $H$ one can induce a new $H$-Hopf module structure as follows.
\begin{proposition}\label{prop.ind}
Let $H$ be a Hopf quasigroup and let $M$ be a right $H$-Hopf module with action $\rom$ and coaction $\roM$. Then $M$ is a right $H$-Hopf module with action
\begin{equation}\label{ind}
\lam =  \rom\circ (\rom\ot H)\circ (M\ot S\ot \mu)\circ (\roM\ot H\ot H)\circ (\roM\ot H),
\end{equation}
and coaction $\roM$.
\end{proposition}
\begin{proof}
In terms of Sweedler's notation, the map $\lam$ reads explicitly, for all $m\in M$ and $h\in H$,
$
\lam(m\ot h) = \left(m\sw 0\act Sm\sw 1\right)\act \left(m\sw 2h\right).
$
The second of equations \eqref{mod1} and the counitality of $\roM$ immediately imply that $\lam(m\ot 1) = m$.
 Write $m\ra h$ for $\lam(m\ot h)$. Then, for all $g,h\in H$, $m\in M$,
$$
(m\ra g)\ra h = \left(\left(\left( m\sw 0\act Sm\sw 1\right)\act \left(m\sw 2g\right)\!\sw 1\right)\act S\left(m\sw 2g\right)\!\sw 2\right) \act \left(\left(m\sw 2g\right)\!\sw 3h\right) ,
$$
by equation \eqref{rom2} in Lemma~\ref{lem.rom1}. Thus the first of equations \eqref{mod1} implies that
 \begin{equation}\label{laM1}
 (m\ra g)\ra h = \left(m\sw 0\act Sm\sw 1\right)\act \left(\left(m\sw 2g\right)h\right).
 \end{equation}
 Therefore,
 $$
 (m\ra Sh\sw 1)\ra h\sw 2 =  \left(m\sw 0\act Sm\sw 1\right)\act \left(\left(m\sw 2Sh\sw 1\right)h\sw 2\right)
  = (m\sw 0\act Sm\sw 1)\act m\sw 2 \eps(h) = m\eps(h),
$$
 where the second  equality follows by the second of equations \eqref{quasi2} and the third one by the second of equations \eqref{mod1}.  Similarly, the first equations in  \eqref{quasi2} and \eqref{mod1} imply that $(m\ra h\sw 1)\ra Sh\sw 2 =  m\eps(h)$. Finally, combining equation \eqref{rom2} in Lemma~\ref{lem.rom1} with the multiplicativity of coproduct we obtain
 \begin{eqnarray*}
 \roM(m\ra h) &=& \left (m\sw 0\act Sm\sw 1\right)\act (m\sw 2h)\sw 1 \ot (m\sw 2h)\sw 2 \\ &=&  \left (m\sw 0\act Sm\sw 1\right)\act (m\sw 2h\sw 1) \ot m\sw 3h\sw 2 
 = m\sw 0\ra h\sw 1\ot m\sw 1 h\sw 2.
 \end{eqnarray*}
 Therefore the compatibility condition \eqref{mod2} between $\lam$ and $\roM$ holds, and $M$ is a right $H$-Hopf module  as stated. 
 \end{proof}
 
One could attempt to iterate the induction procedure described in Proposition~\ref{prop.ind}, i.e.\ apply it again to the $H$-Hopf module with action $\lam$. The antipode axiom and equation \eqref{rom1} in Lemma~\ref{lem.rom1} immediately imply, however, that the iteration must terminate  after the first step, that is $\widehat{\widehat{\varrho}}_M = \lam$. If, on the other hand, the action $\rom$ is associative, then the first of equations \eqref{quasi1} implies that already $\lam = \rom$. 

 \begin{definition}\label{def.qlin}
Let $H$ be a Hopf quasigroup, and let $M$, $N$ be right $H$-Hopf modules. 
A  $\k$-linear map $f: M\to N$ is said to be {\em $H$-quasilinear} provided, it is $H$-linear with respect to the induced actions $\lam$ and $\lan$, defined in Proposition~\ref{prop.ind}. That is, the following diagram
$$
\xymatrix{M\ot H \ar[rr]^-{f\ot H} \ar[d]_{\lam} && N\ot H\ar[d]^{\lan} \\
M \ar[rr]^-f && N}
$$
is commutative. A {\em map of right $H$-Hopf modules} is a map that is both right $H$-colinear and right $H$-quasilinear. The  collection of all right $H$-Hopf modules with their maps forms a category 
which is denoted by $\cM_{\uH}^H$.
\end{definition}

Explicitly and exploring the $H$-colinearity, an $H$-colinear map $f:M\to N$ is $H$-quasilinear (i.e.\ a morphism of $H$-Hopf modules) if and only if,
for all $m\in M$ and $h\in H$,
\begin{equation}\label{qlin}
f\left(\left(m\sw 0 \act Sm\sw 1\right)\act \left(m\sw 2 h\right)\right) = \left( f\left( m\sw 0\right) \act Sm\sw 1\right) \act \left(m\sw 2 h\right).
\end{equation}
 Note that, whenever the $H$-action is associative, the quasilinearity of a map is tantamount with its right $H$-linearity. Thus, if $H$ is a Hopf algebra, then $\cM_{\uH}^H$ is the usual category of right $H$-Hopf modules.

\begin{lemma}\label{lem.rom2}
Let $H$ be a Hopf quasigroup and let $M$ be a right $H$-Hopf module. Then
$M\ot H$ is a right $H$-Hopf module with coaction $M\ot \Delta$ and the diagonal action, i.e., for all $m\in M$, $g,h\in H$,
$$
(m\ot h)\act g := m\act g\sw 1 \ot hg\sw 2.
$$
Furthermore, $\roM$ is a morphism of right $H$-Hopf modules.
\end{lemma}
\begin{proof}
First we need to prove conditions \eqref{mod1} for $M\ot H$. Take any $m\in M$ and $g,h\in H$, and compute
$$
((m\ot h)\act g\sw 1)\act Sg\sw 2 = (m\act g\sw 1)\act Sg\sw 4 \ot (hg\sw 2)Sg\sw 3 
= (m\act g\sw 1)\act Sg\sw 2 \ot h = m\ot h\eps(g),
$$
where the second equality follows by the first of equations \eqref{quasi2}, and the third equality is a consequence of the first of equations  \eqref{mod1}. This proves the first of condition \eqref{mod1} for $M\ot H$. The second one is proven in a similar way, using the second of \eqref{quasi2} and  \eqref{mod1}.

The compatibility condition \eqref{mod2} for $M\ot H$ follows by the compatibility for $M$ and by the multiplicativity of the coproduct (as in the standard Hopf algebra case).

By the coassociativity, the map $\roM$ is right $H$-colinear, so only condition \eqref{qlin} need be checked. Equation \eqref{rom2} implies, for all $m\in M$, $h\in H$,
\begin{eqnarray*}
\roM \left(\left( m\sw 0\act Sm\sw 1\right)\act \left(m\sw 2h\right)\right) &=&  \left( m\sw 0\act Sm\sw 1\right)\act \left(m\sw 2h\right)\!\sw 1\ot \left(m\sw 2h\right)\!\sw 2\\
& =& \left( m\sw 0\act Sm\sw 1\ot 1 \right)\act \left(m\sw 2h\right).
\end{eqnarray*}
On the other hand, the antipode property of $S$ yields,
\begin{eqnarray*}
 \left(\roM \left( m\sw 0\right)\act Sm\sw 1\right)\act \left(m\sw 2h\right) &=& 
 \left(m\sw 0 \act Sm\sw 3 \ot m\sw 1Sm\sw 2\right)\act \left(m\sw 2h\right) \\
 &=&
 \left( m\sw 0\act Sm\sw 1\ot 1 \right)\act \left(m\sw 2h\right).
\end{eqnarray*}
 Therefore, the map $\roM$  is both   right $H$-colinear and $H$-quasilinear as required.
\end{proof}

The category ${}^H_{\uH}\cM$ of left $H$-Hopf modules is defined symmetrically, and the left $H$-Hopf module versions of Lemmas~\ref{lem.rom1} and \ref{lem.rom2}, and of Proposition~\ref{prop.ind} can be similarly proven.

The following theorem contains the main result of this section.

\begin{theorem}\label{thm.qmain}
For any Hopf quasigroup $H$, the categories of left and right $H$-Hopf modules are equivalent to the category of vector spaces.
\end{theorem} 
\begin{proof}
We concentrate on the right $H$-Hopf modules, as the left $H$-Hopf modules can be treated symmetrically. 

Define the {\em induction} functor $F:\vect_\k \to \cM^H_{\uH}$ as follows. For any vector space $V$, $F(V) = V \ot H$ with the $H$-action and coaction defined by, for all $g,h\in H$ and $v\in V$,
$$
(v\ot h)\act g = v\ot hg, \qquad v\ot h\mapsto v\ot h\sw 1\ot h\sw 2.
$$
That is the right $H$-action is $V\ot \mu$ and the coaction is $V\ot \Delta$. Then conditions \eqref{quasi2} ensure that $F(V)$ obeys equations \eqref{mod1}. The multiplicativity of the coproduct yields the satisfaction of condition \eqref{mod2} for $F(V)$. Thus $F(V)$ is a right $H$-Hopf module.
Given a $\k$-linear map $f: V\to W$, define the $H$-colinear and $H$-quasilinear map $F(f)$ by
$$
F(f) : V\ot H\to W\ot H, \qquad v\ot h\mapsto f(v)\ot h.
$$

Define the {\em coinvariants} functor $G:  \cM^H_{\uH}\to \vect_\k$ as follows. For any right $H$-Hopf module $M$,
$$
G(M) = M^{coH} := \{ m\in M\; |\; \roM(m) = m\ot 1\}.
$$
On morphisms $G$ treats an $H$-colinear and $H$-quasilinear map $f: M\to N$ as the underlying linear transformation and restricts it to $M^{coH}$.  The image of the restricted map is a subspace of $N^{coH}$ since $f$ is  $H$-colinear. 
\begin{lemma}\label{lem.adj}
The induction functor $F$ is left adjoint to the coinvariants functor $G$.
\end{lemma}
\begin{proof}
We need to construct the unit and counit of adjunction. For any vector space $V$, define
$$
\eta_V : V \to (V\ot H)^{coH}, \qquad v\mapsto v\ot 1.
$$
This is obviously a $\k$-linear map and it is natural in $V$, thus the family of all such $\eta_V$ forms a natural transformation $\eta: \id_{\vect_\k} \to G\circ F$. 

For any right $H$-Hopf module define,
$$
\sigma_M : M^{coH}\ot H \to M, \qquad m\ot h\mapsto m\act h.
$$
The $H$-colinearity of $\sigma_M$  follows by the compatibility condition \eqref{mod2}, as, for all $m\in M^{coH}$ and $h\in H$,
$$
\roM(\sigma_M(m\ot h)) = (m\act h)\sw 0 \ot (m\act h)\sw 1 = m\act h\sw 1\ot h\sw 2 = \sigma_M( m \ot h\sw 1)\ot h\sw 2.
$$
 Next take any $g,h\in H$ and $m\in M^{coH}$ and compute
\begin{eqnarray*}
\sigma_M\left(\left(\left( m\ot h\sw 1\right)\act Sh\sw 2\right) \act \left( h\sw 3g\right)\right) &=& \sigma_M \left(m \ot  \left(h\sw 1 Sh\sw 2\right) \left(h\sw 3g\right)\right) \\
&=& \sigma_M(m\ot hg) = m\act (hg),
\end{eqnarray*}
where the definition of the antipode have been used to derive the third equality. On the other hand,
$$
\left(\sigma_M\left( m\ot h\sw 1\right)\act Sh\sw 2\right) \act \left( h\sw 3g\right) =\left(\left( m\act h\sw 1\right)\act Sh\sw 2\right) \act \left( h\sw 3g\right) = m\act (hg),
$$
by the first of equations \eqref{mod1}. Therefore, $\sigma_M$ is also a right $H$-quasilinear map as required.

Take any $H$-colinear and $H$-quasilinear map $f: M\to N$ and compute, for all $m\in M^{coH}$ and $h\in H$,
\begin{eqnarray*}
f(\sigma_M(m\ot h)) &=& f(m\act h) = f\left(\left(m\sw 0\act Sm\sw 1\right)\act (m\sw 2 h)\right)\\
&=& \left( f\left(m\sw 0\right)\act Sm\sw 1\right)\act (m\sw 2 h) = f(m)\act h = \sigma_N(f(m)\ot h),
\end{eqnarray*}
where the second and fourth equalities follow by the fact that $\roM(m) = m\ot 1$ and by the unitality of actions and products, while the third equality is a consequence of $H$-quasilinearity of $f$. This proves that $\sigma_M$ is natural in $M$, hence the family of all the $\sigma_M$ forms  a natural transformation $\sigma: F\circ G\to \id_{\cM^H_{\uH}}$. Finally, the triangular identities for the unit and counit of adjunction follow by the unitality of both the $H$-action  and product in $H$. 
\end{proof}

In view of Lemma~\ref{lem.adj} we need to prove that the unit and counit of adjunction constructed in its proof are natural isomorphisms. For any vector space $V$, the inverse of $\eta_V$ is given by
$$
\eta^{-1}_V: (V\ot H)^{coH} \to V, \qquad \sum_i v^i\ot h^i\mapsto \sum_i v^i\eps(h^i).
$$
The map $\eta^{-1}_V$ is clearly natural in $V$, and that the composite $\eta^{-1}_V\circ \eta_V$ is the identity linear transformation on $V$ follows immediately by $\eps(1)=1$. To compute the other composite, first observe that $\sum_i v^i\ot h^i \in  (V\ot H)^{coH}$ if and only if
$$
\sum_i v^i\ot h^i\sw 1\ot h^i\sw 2 = \sum_i v^i\ot h^i\ot 1.
$$
Applying the counit to the middle factor we thus obtain that the elements of $ (V\ot H)^{coH}$ are characterised by the equation $\sum_i v^i\ot h^i = \sum_i v^i\eps(h^i)\ot 1$. Using this characterisation one immediately finds that the composite $\eta_V\circ \eta^{-1}_V$ is the identity on $(V\ot H)^{coH}$. All this is the same as in the standard Hopf algebra case.

The construction of the natural inverse of $\sigma$ is slightly more involved. We claim that, for all $H$-Hopf modules $M$, the inverse of $\sigma_M$ is given by
$$
\sigma^{-1}_M: M\to M^{coH}\ot H, \qquad m\mapsto m\sw 0\act Sm\sw 1\ot m\sw 2.
$$
By equation \eqref{rom1} in Lemma~\ref{lem.rom1},  the map $\sigma_M^{-1}(M)\subseteq M^{coH}\ot H$, hence $\sigma^{-1}_M$ is well-defined (as a $\k$-linear map). By the coassociativity of the coproduct,  $\sigma^{-1}_M$ is a right $H$-colinear map.
In view of equation \eqref{rom2} in Lemma~\ref{lem.rom1}, for all $h\in H$ and $m\in M$,
\begin{eqnarray*}
\sigma_M^{-1} \!\!\!\!\!\!\!\!&&\!\! \left(\left(m\sw 0 \act Sm\sw 1\right)\act \left(m\sw 2 h\right)\right)\\
&& =  \left(\left(m\sw 0 \act Sm\sw 1\right) \act \left(m\sw 2 h\right)\!\sw 1 \right)\act S\left(m\sw 2 h\right)\!\sw 2 \ot \left(m\sw 2h\right)\!\sw 3
= m\sw 0\act Sm\sw 1 \ot m\sw 2h,
\end{eqnarray*}
where the last equality follows by the first of equations \eqref{mod1}. On the other hand,
\begin{eqnarray*}
\left(\sigma^{-1}_M\left( m\sw 0\right) \act Sm\sw 1\right) \act \left(m\sw 2 h\right) &=& \left(\left( m\sw 0\act Sm\sw 1\ot m\sw 2\right)\act Sm\sw 3\right)\act \left(m\sw 4h\right)\\
&=& m\sw 0\act Sm\sw 1\ot \left(m\sw 2 Sm\sw 3\right) \left(m\sw 4h\right)\\
&=& m\sw 0\act Sm\sw 1 \ot m\sw 2h,
\end{eqnarray*}
where the antipode property was used to obtain the final equality. Therefore, $\sigma^{-1}_M$ is also an $H$-quasilinear map. 

Take any $H$-colinear and $H$-quasilinear map $f:M\to N$. Since $f$ is $H$-colinear, 
$$
\sigma^{-1}_N\circ f(m) = f(m)\sw 0 \act Sf(m)\sw 1 \ot f(m)\sw 2 = f(m\sw 0) \act Sm\sw 1 \ot m\sw 2,
$$
 for all $m\in M$.  On the other hand, by the antipode property, $H$-quasilinearity of $f$, and the antipode property again,
\begin{eqnarray*}
 f\left(m\sw 0\act Sm\sw 1\right) &=& f\left(\left(m\sw 0\act Sm\sw 1\right)\act \left(m\sw 2 Sm\sw 3\right)\right)\\
 & =& \left(f(m\sw 0)\act Sm\sw 1\right)\act \left(m\sw 2 Sm\sw 3\right) = f(m\sw 0)\act Sm\sw 1.
 \end{eqnarray*}
 Therefore,
 $$
 (f \ot H)\circ \sigma^{-1}_M (m) = f(m\sw 0\act Sm\sw 1) \ot m\sw 2 = f(m\sw 0)\act Sm\sw 1\ot m\sw 2 = \sigma^{-1}_N\circ f(m),
 $$
i.e.\ the map $\sigma^{-1}_M$ is natural in $M$. 

Finally we need to prove that, for every $H$-Hopf module $M$,  $\sigma^{-1}_M$ is the inverse of $\sigma_M$. For any $m\in M$,
$$
\sigma_M(\sigma_M^{-1}(m)) = (m\sw 0\act Sm\sw 1)\act m\sw 2 = m,
$$
by the second of equalities \eqref{mod1} and the counitality of the coaction. For all $m\in M^{coH}$ and $h\in H$,
$$
\sigma^{-1}_M(\sigma_M(m\ot h)) = (m\act h)\sw 0\act S(m\act h)\sw 1 \ot (m\act h)\sw 2
= (m\act h\sw 1)\act Sh\sw 2\ot h\sw 3 = m\ot h,
$$
where the second equality is obtained by combining equation \eqref{mod2} with the fact that $\roM(m) = m\ot 1$, and the third equality follows by the first of equations \eqref{mod1} and by the counitality of the coproduct. 

This completes the proof of the category equivalence $\cM^H_{\uH}\simeq \vect_\k$. The  equivalence of vector spaces to left $H$-Hopf modules follows by symmetric arguments.
\end{proof}

Directly from the proof of Theorem~\ref{thm.qmain} one draws the following 
\begin{corollary}\label{cor.qmain}
Let $H$ be a Hopf quasigroup, and view $H\ot H$ as a right $H$-Hopf module as in Lemma~\ref{lem.rom2}. Let $F: \vect_\k \to \cM^H_{\uH}$ be the induction functor described in the proof of Theorem~\ref{thm.qmain}. Then the right Galois map $\beta$ given in Definition~\ref{def.Gal} is an isomorphism of $H$-Hopf modules $F(H) \to H\ot H$.
\end{corollary}
\begin{proof}
Observe that $(H\ot H)^{coH} = H\ot 1$. Then $\beta = \sigma_{H\ot H}$, where $\sigma$ is the counit of the induction-coinvariant adjunction described in the proof of Theorem~\ref{thm.qmain}. Thus $\beta$ is an isomorphism of $H$-Hopf modules.
\end{proof}

In a symmetric way the left Galois map $\gamma$ is an isomorphism of left $H$-Hopf modules.

\subsection{Hopf modules over Hopf coquasigroups}
Considerations regarding Hopf modules over Hopf coquasigroups are dual to these regarding Hopf quasigroups, so we merely outline main observations without giving detailed proofs.

\begin{definition}\label{def.comod}
Let $H$ be a Hopf coquasigroup. Let $M$ be a counital right $H$-comodule with coaction $\roM$ and an associative unital right $H$-module with action $\varrho_M$. $M$ is said to be a {\em right $H$-Hopf module} if
\begin{equation}\label{comod1}
(M\ot \mu)\circ (M\ot H\ot S)\circ (\roM \ot H)\circ \roM   = M\ot 1 =  (M\ot \mu)\circ (M\ot S\ot H)\circ(\roM \ot H)\circ \roM
\end{equation}
and
\begin{equation}\label{comod2}
\roM\circ \rom = (\rom \ot \mu)\circ (M\ot \flip \ot H)\circ (\roM\ot \Delta).
\end{equation}
\end{definition}
The conditions \eqref{comod1} and \eqref{comod2} are obtained by the dualisation of conditions \eqref{mod1} and \eqref{mod2}. Note that \eqref{mod2} is self-dual hence conditions \eqref{mod2} and \eqref{comod2} coincide. 
By the same dualisation procedure applied 
 to Proposition~\ref{prop.ind} one obtains
\begin{proposition}\label{prop.coind}
Let $H$ be a Hopf coquasigroup, and let $M$ be a right $H$-Hopf modules with action $\rom$ and coaction $\roM$. Then $M$ is a right $H$-Hopf module with action $\rom$ and coaction
$$
\laM = (\rom\ot H)\circ (\rom\ot H\ot H)\circ (M\ot S\ot \Delta)\circ (\roM\ot H)\circ \roM.
$$
\end{proposition}

As in the case of modules over Hopf coquasigroups, the iteration of coactions by the procedure described in Proposition~\ref{prop.coind} terminates after the first step, i.e. $\widehat{\widehat{\varrho}}^M = \laM$. 

\begin{definition}\label{def.coqlin}
Let $H$ be a Hopf coquasigroup, and let $M$, $N$ be right $H$-Hopf modules. 
A  $\k$-linear map $f: M\to N$ is said to be {\em right $H$-quasicolinear} provided it is right $H$-colinear with respect to the induced coactions $\laM$, $\laN$; see Proposition~\ref{prop.coind}.
The  collection of all right $H$-Hopf modules  with right $H$-linear and $H$-quasicolinear maps as morphisms forms a category 
which is denoted by $\cM^{\oH}_H$.
\end{definition}

The category of left $H$-Hopf modules
is defined symmetrically.
Dualising the statements and proofs of Lemmas~\ref{lem.rom1} and \ref{lem.rom2} one obtains the following
\begin{lemma}\label{lem.corom}
Let $H$ be a Hopf coquasigroup and let $M$ be a right $H$-Hopf module. 
\begin{zlist}
\item For all $m\in M$ and $h\in H$,
$$
(m\act h)\sw 0\act S(m\act h)\sw 1 = m\sw 0\act Sm\sw 1\eps(h).
$$
\item
 $M\ot H$ is a right $H$-Hopf module with action $M\ot \mu$ and the diagonal coaction, i.e.\ for all $m\in M$, $h\in H$,
$$
m\ot h \mapsto  m\sw 0\ot  h\sw 1 \ot m\sw 1h\sw 2.
$$
Furthermore, $\rom$ is a morphism of $H$-Hopf modules.
\end{zlist}
\end{lemma}

The following is the Hopf coquasigroup version of Theorem~\ref{thm.qmain}.
\begin{theorem}\label{thm.coqmain}
For any Hopf coquasigroup $H$, the categories of left and right $H$-Hopf modules are equivalent to the category of vector spaces.
\end{theorem} 
\begin{proof}
This is dual to Theorem~\ref{thm.qmain}. In the case of right $H$-Hopf modules, the equivalence is given by the induction functor $F= -\ot H:\vect_\k \to \cM_H^{\oH}$ defined by the same formulae as in the proof of Theorem~\ref{thm.qmain}. Its inverse equivalence is the {\em $H$-invariants} functor $G: \cM_H^{\oH}\to \vect_\k$ defined on objects as $G(M) = M^H$, where $M^H$ is given by the coequaliser
$$
\xymatrix{ M\otimes H
\ar@<0.5ex>[rr]^-{\rom}\ar@<-0.5ex>[rr]_-{M\ot \eps} & &
M \ar[rr]^-{\pi_M} &&M^H .}
$$
On morphisms $f: M\to N$, $G(f)$ is defined by the formula $G(f)\circ \pi_M = \pi_N\circ f$ (the existence and the uniqueness of $G(f)$ follow by the fact that $f$  is right $H$-linear and by the universal property of coequalisers). $G$ is the left adjoint of $F$. The unit of adjunction is defined by, for all $H$-Hopf modules $M$,
$$
\eta_M : M\to M^H\ot H,\qquad m\mapsto \pi_M(m\sw 0)\ot m\sw 1.
$$
The counit of adjunction is, for all vector spaces $V$,
$$
\sigma_V: (V\ot H)^H \to V, \qquad \pi_{V\ot H}(v\ot h)\mapsto v\eps(h),
$$ 
which is well-defined by the multiplicativity of the counit $\eps$ and by the universal property of coequalisers. The inverse of $\sigma_V$ is
$$
\sigma_V^{-1}: V\to (V\ot H)^H, \qquad v\mapsto \pi_{V\ot H}(v\ot 1),
$$
and the inverse of $\eta_M$ is
$$
\eta^{-1}_M : M^H\ot H\to M, \qquad \pi_M(m)\ot h\mapsto m\sw 0\act (Sm\sw 1)h,
$$ 
which is well-defined by the antimultiplicativity of $S$ and by its antipode property.
\end{proof}

Dually to Corollary~\ref{cor.qmain}, one proves that the left Galois map $\gamma$ given in Definition~\ref{def.Gal} is an isomorphism of right $H$-Hopf modules $H\ot H\to F(H)$, where $H\ot H$ is a right $H$-Hopf module as in Lemma~\ref{lem.corom}(2).

Note that if $H$ is a Hopf algebra, then the $H$-coinvariants and $H$-invariants functors from the category of  right $H$-Hopf modules to $\vect_\k$ are naturally isomorphic. For a right $H$-Hopf module $M$ the isomorphism $M^{co H}\to M^H$ and its inverse are, for all $n\in M^{co H}$ and $m\in M$,
$$
n\mapsto \pi(n), \qquad \pi(m)\mapsto m\sw 0\act Sm\sw 1.
$$
Under this isomorpism the unit (resp.\ counit) of adjunction constructed in the proof of Theorem~\ref{thm.coqmain} coincides with the inverse of the counit (resp.\ unit) of adjunction describd in the proof of Theorem~\ref{thm.qmain}.


\begin{thebibliography}{99}
\bibitem{BohNil:wea}  G.\ B\"ohm,  F.\ Nill and K.\  Szlach\'anyi,
        {\em Weak  Hopf algebras I. Integral theory and $C^*$-structure},
         J.\ Algebra 221, 385--438 (1999).

\bibitem{Dri:qua} V.G.\ Drinfeld, {\em Quasi-Hopf algebras}, Leningrad Math.\ J.\ 1, 1419--1457 (1990).

\bibitem{KliMaj:Hop} J.\ Klim and S.\ Majid, {\em Hopf quasigroups and the algebraic 7-sphere}, Preprint arXiv:0906.5026v2 (2009). 

\bibitem{Swe:Hop} M.E.\ Sweedler, {\em Hopf Algebras}, W.A.\ Benjamin, New York, (1969).


\bibitem{Van:mul} A.\ Van Daele, {\em Multiplier Hopf algebras}, Trans.\ Amer.\ Math.\ Soc.\ 342,  917--932 (1994).

\end{thebibliography}
 \end{document}